\documentclass[12pt,a4paper,draft]{amsart}
\usepackage[english]{babel}
\usepackage{a4wide,amssymb,xypic}
\newtheorem{thm}{Theorem}[section]
\newtheorem{corol}[thm]{Corollary}
\newtheorem{lemma}[thm]{Lemma}
\newtheorem{prop}[thm]{Proposition}
\newtheorem{defin}[thm]{Definition}

\theoremstyle{remark}
\newtheorem{rema}[thm]{Remark}
 \newenvironment{remark}{\begin{rema}}{\hfill\hspace{1pt}$\triangle$\end{rema}}
\newtheorem{exe}[thm]{Example}
 \newenvironment{example}{\begin{exe}}{\hfill\hspace{1pt}$\triangle$\end{exe}}

\newcommand{\Oc}{\mathcal O}

\newcommand{\PP}{{\mathbb P}}

\newcommand{\D}{\mathcal{D}}

\newcommand{\ccR}{\mathcal{R}}
\newcommand{\R}{\mathbb R}
\newcommand{\C}{\mathbb C}
\newcommand{\Q}{\mathbb Q}

  \def\K{{\mathcal K}}

\def\Id{\operatorname{Id}}
\def\rk{\operatorname{rk}}
\def\vol{\operatorname{vol}}

\def\End{\operatorname{End}}
\def\dim{\operatorname{dim}}

\newcommand\LHom{\mbox{\it Hom}}

\newcommand\grass{\operatorname{Gr}}
\newcommand\hgrass{{\mathfrak{Gr}}}
\newcommand{\cO}{{\mathcal O}}

\newcommand{\fE}{{\mathfrak E}}
\newcommand{\fF}{{\mathfrak F}}
\newcommand{\fG}{{\mathfrak G}}
\newcommand{\fQ}{{\mathfrak Q}}
\newcommand{\fS}{{\mathfrak S}}

\newcommand{\fL}{{\mathfrak L}}
\newcommand{\fK}{{\mathfrak K}}

\newcommand{\OPQ}[1]{{\mathcal O_{\PP Q_{#1}}}}

\begin{document}
\rightline{SISSA Preprint 17/2006/fm}
\rightline{ESI Preprint 1847}
\rightline{\tt arXiv:math.DG/0605659}
\vfill
\title[Metrics on semistable and numerically effective Higgs bundles]{Metrics on semistable and numerically \\[8pt] effective Higgs bundles}
\date{Revised 19 September 2006}
\thanks{This work has been
partially supported by the European  Union through the FP6 Marie
Curie Research and Training Network {\sc enigma} (Contract number
MRTN-CT-2004-5652), by the Italian National Project ``Geometria delle variet\`a algebriche,'' by the Spanish {\scriptsize
DGES} through the research project BFM2003-00097 and by ``Junta de Castilla y Le{\'o}n''
through the research project SA114/04.
 }
 \maketitle \thispagestyle{empty}
\begin{center}{\sc Ugo Bruzzo} \\
Scuola Internazionale Superiore di Studi Avanzati,\\ Via Beirut
2-4, 34013 Trieste, Italia;\\[2pt]
Istituto Nazionale di Fisica Nucleare, Sezione di Trieste \\[2pt]  E-mail: {\tt bruzzo@sissa.it} \\[10pt]
{\sc Beatriz Gra\~na Otero} \\
Departamento de Matem\'aticas and Instituto de F\'\i sica \\ Fundamental y Matem\'aticas, Universidad de Salamanca,
\\ Plaza de la Merced 1-4, 37008 Salamanca, Espa\~na\\ E-mail: {\tt beagra@usal.es}
\end{center}

\begin{abstract} We provide notions of numerical effectiveness
and numerical flatness for Higgs vector bundles on compact
K\"ahler manifolds in terms of   fibre metrics. We prove several
properties of bundles satisfying such conditions and in particular
we show that numerically flat Higgs bundles have vanishing Chern
classes, and that they admit filtrations whose quotients are
stable flat Higgs bundles. We compare these definitions with those
previously given in the case of projective varieties.
Finally we  study the relations between
numerically effectiveness and semistability, establishing
semistability criteria for Higgs bundles on projective manifolds of
any dimension.
\end{abstract}
\vfill
\newpage

\section{Introduction}
Numerical effectiveness is a natural generalization of the notion
of ampleness: a line bundle $L$ on a projective variety $X$ is
said to be numerical effective (nef) if $c_1(L) \cdot [C] \geq 0$
for every irreducible curve $C$ in $X$. A vector bundle $E$ is
said to be nef if the hyperplane bundle $\Oc_{\PP E}(1)$ on the
projectivized bundle $\PP E$ is nef. A definition of numerical
effectiveness can also be given in the case of Hermitian or
K\"ahler manifolds in terms of (possibly singular) Hermitian fibre
metrics.

The classification of the complex projective varieties or compact
K\"ahler manifolds whose tangent bundle is numerically effective
\cite{DPS94,CP91} yields a far-reaching generalization of the
Hartshorne-Frankel-Mori theorem, the latter stating that  the complex
projective $n$-space is the only projective $n$-variety whose
tangent bundle is ample. Manifolds whose cotangent bundle is nef
were studied by Kratz \cite{K97}.

The concept of nefness may also be used to provide a
characterization of semistable bundles. Several criteria of this
type have appeared recently \cite{BHR05,BB04,BSch05}, all
generalizing Miyaoka's result according to which a vector bundle
$E$ on a smooth projective curve $X$ is semistable if and only if
the numerical class $\lambda = c_1(\Oc_{\PP E}(1)) - \frac{1}{r}
\pi^{\ast}(c_1(E))$ is nef, where $\pi: \PP E \to X$ is the
projection. Also results by Gieseker \cite{G79}, generalized in
\cite{BHR05,BG05}, relate the notions of nefness and
semistability.

A notion of nefness for Higgs bundles was introduced \cite{BG05}
in the case of projective varieties, and the basic properties of
such bundles were there studied. In particular, it was proved there that
numerically flat Higgs bundles (i.e., Higgs bundles that are
numerically flat together with their duals) have vanishing Chern
classes. One could stress here a similarity with
the Bogomolov inequality: Higgs bundles that are semistable as
Higgs bundles but not as ordinary bundles nevertheless satisfy
Bogomolov's inequality.

The main purpose of this paper is to extend these constructions to
the case when the base variety is a compact K\"ahler manifold.
The guiding idea here is to give a definition of nefness in terms
of fibre metrics on the bundles, following
Demailly-Peternell-Schneider and de Cataldo \cite{DPS94,dC98}. One
takes the Higgs structure  into account  by replacing the Chern
connection associated with the fibre metrics by a connection
(introduced by Hitchin \cite{Hitch} and later extensively used by Simpson \cite{S92}) whose definition involves the Higgs field.

We now describe the contents of this paper. In Section \ref{2},
the Hitchin-Kobayashi correspondence for Higgs bundles is extended
to a relation between semistability and the existence of
approximate Hermitian-Yang-Mills structures. This will be used in
the ensuing sections for our study of numerically effectiveness of
Higgs bundles. We also give Gauss-Codazzi equations associated
with extensions of Higgs bundles.

In Section \ref{3} we give our definition of numerically
effectiveness and prove several related basic properties. The main
result of this section is a characterization of numerically flat
Higgs bundles on K\"ahler manifolds as those Higgs bundles which
admit a filtration whose quotients are flat stable Higgs bundles
(flat in a sense to be specified later). This implies that the
Chern classes of numerically flat Higgs bundles vanish. We also
establish some relations between semistability and numerical
effectiveness.

In Section \ref{4} we compare these results with our previous work
in the case of projective varieties. This will also allow us to
complete some of the results given in \cite{BG05}. In particular, we
provide new criteria for the
characterization (in terms of the
numerical effectiveness of certain associated  Higgs bundles) of those Higgs
bundles on projective manifolds that are semistable and
satisfy the condition $\Delta(E)=c_2(E)-\frac{r-1}{2r}c_1(E)^2=0$.

We give now the basic definitions concerning Higgs bundles. Let
$X$ be an $n$-dimensional compact K\"ahler manifold, with K\"ahler
form $\omega$. Given a coherent sheaf $F$ on $X$, we denote by
$\deg(F)$ its \emph{degree} $$\deg(F) = \int_X c_1(F)\cdot \omega^{n-1}$$
and if $r=\rk(F)>0$ we introduce its \emph{slope}
$$\mu(F) = \frac{\deg(F)}{r}.$$

\begin{defin} A Higgs sheaf $\fE$ on $X$ is a coherent sheaf $E$ on $X$
endowed with a morphism $\phi \colon E \to E \otimes \Omega_X$ of
$\Oc_X$-modules such that the morphism $\phi\wedge\phi : E \to E
\otimes \Omega^2_X$ vanishes, where $\Omega_X$ is the cotangent
sheaf to $X$. A Higgs subsheaf $F$ of a Higgs sheaf $\fE=(E,\phi)$
is a subsheaf of $E$ such that $\phi(F)\subset F\otimes\Omega_X$.
A Higgs bundle is a Higgs sheaf $\fE $ such that $E$ is a
locally-free $\Oc_X$-module.
\end{defin}

There exists a stability condition for Higgs sheaves, analogous to
that for ordinary sheaves, which makes reference only to
$\phi$-invariant subsheaves.

\begin{defin}  A Higgs sheaf $\fE=(E,\phi) $ on $X$ is semistable
(resp.~stable) if $E$ is torsion-free, and $\mu(F)\le \mu(E)$
(resp. $\mu(F)< \mu(E)$) for every proper nontrivial Higgs
subsheaf $F$ of $\fE$.\end{defin}

{\bf Acknowledgements.} The authors thank M.S.~Narasimhan for
drawing de Cataldo's paper \cite{dC98} to their attention, and the referee
for helping to improve the exposition. This
paper was mostly done during a stay of both authors at the Tata
Institute for Fundamental Research in Mumbai, India; they
express their warm thanks for the hospitality and support. The
stay was also made possible by grants from Istituto Nazionale di
Alta Matematica and Universidad de Salamanca. The paper was finalized
while both authors were visiting the Erwin Schr\"odinger Institut f\"ur Mathematische
Physik in Vienna.

\section{Metrics and connections on semistable Higgs
bundles}\label{2} This section mostly deals with an application to
Higgs bundles of the beautiful ideas underlying the so-called
Hitchin-Kobayashi correspondence. Kobayashi and L\"ubke
\cite{koba,Lub83} proved that the sheaf
of sections of a holomorphic Hermitian (ordinary) vector bundle
satisfying the Hermitian-Yang-Mills condition is polystable (i.e.,
it is a direct sum of stable sheaves having the same slope)
A  converse result was proved first by Donaldson
in the projective case \cite{Don85,Don87}, and then by Uhlenbeck and Yau
in the compact K\"ahler case \cite{UY}; they showed that a stable bundle admits
a (unique up to homotheties) Hermitian metric which satisfies the
Hermitian-Yang-Mills condition. One can also show that if
Hermitian bundle satisfies the Hermitian-Yang-Mills condition in
an approximate sense, then it is semistable,  while the converse may
be proved if $X$ is projective \cite{koba},

Later Simpson \cite{S88,S92} proved a Hitchin-Kobayashi
correspondence for Higgs bundles. Given a Higgs bundle equipped
with an Hermitian metric, one defines a natural connection (that
we call the \emph{Hitchin-Simpson connection}); when this satisfies the Hermitian-Yang-Mills condition then the
Higgs bundle is polystable, and \emph{vice versa}
(this connection was originally introduced by Hitchin in \cite{Hitch}). Our main aim in
this section is to show that whenever an Hermitian Higgs bundle
satisfies an approximate Hermitian-Yang-Mills condition, then it
is semistable.  To this end we shall need to prove a vanishing result.

\subsection{Main definitions}
Let $\fE = (E, \phi)$ be a Higgs bundle over a complex manifold
$X$ equipped with an Hermitian  fibre metric $h$. Then there is on
$E$ a unique connection $D_{(E,h)}$ which is compatible with both
the metric $h$ and the holomorphic structure of $E$. This is often
called the \emph{Chern connection} of the Hermitian bundle
$(E,h)$.

Now let  $\bar\phi$ be the adjoint of the morphism $\phi$ with
respect to the metric $h$, i.e., the morphism $\bar\phi: E
\rightarrow \Omega^{0,1}\otimes E $ such that $$h(s, \phi(t)) =
h(\bar\phi(s), t)$$ for all sections $s,t$ of $E$. It is easy to
check that the operator
\begin{equation}\label{e:simpsonconnection}
\D_{(\fE,h)} = D_{(E,h)} + \phi + \bar\phi
\end{equation}
defines a connection on the bundle $E$. One should notice that
this connection is neither compatible with the holomorphic
structure of $E$, nor with the Hermitian metric $h$.
\begin{defin} \label{simpsonconnection}
The connection \eqref{e:simpsonconnection} is called the
\emph{Hitchin-Simpson connection} of the Hermitian Higgs bundle $(\fE,h)$.
Its curvature will be denoted by $\ccR_{(\fE,h)} = \D_{(\fE,h)}
\circ \D_{(\fE,h)}$. If this curvature vanishes, we say that
$(\fE, h)$ is \emph{Hermitian flat}.
\end{defin}
If $\fE$ is a  Higgs line bundle the notion of Hermitian flatness coincides
with the usual one.

Let us now assume that $X$ has a K\"ahler metric, with K\"ahler
form $\omega$. We shall denote by
$\mathcal{K}_{(\fE,h)}\in\End(E)$ the mean curvature of the
Hitchin-Simpson connection. This is defined as usual: if we consider
wedging by the K\"ahler 2-form as a morphism $\mathcal A^p\to
\mathcal A^{p+2}$ (where $\mathcal A^p$ is the sheaf of $\C$-valued smooth $p$-forms on $X$), and denote by $\Lambda: \mathcal A^p\to
\mathcal A^{p-2}$ its adjoint, then $\mathcal{K}_{(\fE,h)}= i
\Lambda \ccR_{(\fE,h)}$. Locally, if we write
$$ \ccR_{(\fE,h)} = \tfrac12 \sum_{i,j,\alpha,\beta}(\ccR_{(\fE,h)})^{i}_{\ j\alpha\beta}\,
e_i\otimes (e^\ast)^j \otimes dx^\alpha\wedge dx^\beta$$ in terms
of a local basis $\{e_i\}$ of sections of $E$, and local real
coordinates $\{x^1,\dots,x^{2n}\}$ on $X$, we have
$$(\mathcal{K}_{(\fE,h)}) ^{l}_{\ j} =- i \sum_{\alpha,\beta}
(\omega^{-1})^{\alpha\beta}\,(\ccR_{(\fE,h)})^{l}_{\ j\alpha\beta}\,.$$
We shall regard the mean curvature as a bilinear form on sections
of $E$ by letting
$$\mathcal{K}_{(\fE,h)}(s,t)= h(\mathcal{K}_{(\fE,h)}(s),t)\,.$$

Let us recall, for comparison and   further use, the form that the
Hitchin-Kobayashi correspondence acquires for Higgs bundles
\cite[Thm.~1]{S92}.

\begin{thm} A   Higgs vector
bundle $\fE =(E, \phi)$ over  a compact K{\"a}hler manifold  is
polystable if and only if it admits an Hermitian metric $h$ such
that the curvature of  the Hitchin-Simpson connection of $(\fE,h)$
satisfies the  Hermitian-Yang-Mills condition
$$ \mathcal{K}_{(\fE,h)}= c \cdot \Id _E $$
for some constant real number $c$.
\end{thm}
The constant $c$ is related to the slope of $E$:
\begin{equation}\label{c}
c=\frac{2n\pi}{n! \vol(X)} \mu(E)
\end{equation}
where $n=\dim(X)$ and
$$\vol(X) = \frac{1}{n!}\int_X \omega^n\,.$$

As in the case of ordinary bundles, the semistability of a Higgs
bundle may be related to the existence of an approximate
Hermitian-Yang-Mills structure, which is introduced by replacing
the Chern connection by the Hitchin-Simpson connection. Before doing that,
we need to write the Gauss-Codazzi equations associated with an
extension of Hermitian Higgs bundles.

\subsection{Metrics on Higgs subbundles and quotient bundles}
\label{codazzi}
Given a rank $r$ Higgs bundle $\fE=(E,\phi)$, let $\fS =(S,
\phi_{\vert S})$ be a rank $p$ Higgs subbundle (thus, $\phi_{\vert
S}(S)\subset S\otimes \Omega^1$). The quotient bundle $Q$ has an
induced Higgs field $\phi_Q$, so that one has
 a quotient Higgs bundle $\fQ = (Q, \phi_Q)$.
If $h$ is an Hermitian metric on $\fE$ the orthogonal complement $S^\perp$ is isomorphic to $Q$ in the $C^\infty$ category,
including the extra structure as Higgs bundles, so that one gets a
$C^\infty$ decomposition of $\fE$ into Higgs subbundles
$$ \fE \simeq \fS \oplus \fS^{\perp}\,.$$
Note that this endows the bundle $Q$ with an Hermitian metric
$h_Q$.

Let $e_1, \cdots, e_p, e_{p+1}, \cdots, e_r$ be a local
$C^{\infty}$ unitary frame field for $E$,  such that the first $p$
sections form a local unitary frame for $S$, and the last $r-p$
ones yield a local unitary frame for $\fS^{\perp}  \simeq \fQ$. We
shall use the following notation for the indices ranging over the
sections of the various bundles:
$$ 1\leq a,b,c \leq p, \qquad p+1 \leq \lambda, \mu, \nu \leq r,
\qquad 1 \leq i, j, k \leq r. $$
We introduce local connection 1-forms by letting
\begin{eqnarray*}
\D_{(\fE,h)}(e_a) & = & \sum_{b=1}^p \omega^b_a \otimes e_b +
\sum_{\lambda=p+1}^r\omega^{\lambda}_a \otimes e_{\lambda}
\end{eqnarray*}
\begin{eqnarray*}
\D_{(\fE,h)}(e_{\lambda}) & = & \sum_{a=1}^p \omega^a_{\lambda}
\otimes e_a + \sum_{\mu=p+1}^r\omega^{\mu}_{\lambda} \otimes
e_{\mu}\,.
\end{eqnarray*}
Let us furthermore set
\begin{equation}\label{fund} A(e_a) = \sum_{\lambda=p+1}^r\omega^{\lambda}_a \otimes e_{\lambda}
\qquad\mbox{and}\qquad
B(e_{\lambda}) = \sum_{a=1}^{p} \omega^a_{\lambda} \otimes e_a \,.\end{equation}

\begin{prop}
\begin{enumerate}
    \item The 1-forms $\omega^b_a$ are local connection 1-forms
    for the Hitchin-Simpson connection $\D_{(\fS,h_S)}$ of $(\fS,h_S)$;
    \item the map defined by $A$ is a $(1,0)$-form with values in
    $\LHom(S, S^{\perp})$;
    \item   the map $B$ is a $(0,1)$-form with values in
    $\LHom(S^{\perp},S)$.
    \end{enumerate}
\end{prop}
\begin{proof} The first statement follows from the fact that
$\fS$ is a Higgs subbundle of $\fE$. So we have $\D_{(\fE,h)}(e_a)
= \D_{(\fS,h_S)}(e_a) + A(e_a)$. Since both connections satisfy
the Leibniz rule, $A$ is $f$-linear. The same is true for $B$. Moreover, $\D_{(\fE,h)}$ and
$ \D_{(\fS,h_S)}$ have the same $(0,1)$ part, namely,
$\bar\partial+\bar\phi$, so that $A$ is of type (1,0).
Analogously, $B$ is of type $(0,1)$.
\end{proof}
We want now to relate the curvatures of the Hitchin-Simpson connections on
the Higgs bundles $\fS$ and $\fE$ (namely, we want to write the
relevant equations of the Gauss-Codazzi type). Since (for
simplicity we omit here the tensor product symbols and the ranges
of summation over the indices)
\begin{eqnarray*}
D_{(\fE, h)}^2(e_a) & = & D_{(\fE,h)}[\, \Sigma_b \, \omega^b_a \,
e_b + \Sigma_{\lambda} \, \omega^{\lambda}_a \, e_{\lambda}] =  \\
& = & \Sigma_b [d \omega_a^b \, e_b - \omega_a^b \wedge (\Sigma_c
\, \omega_b^c\, e_c + \Sigma_{\mu} \, \omega_b^\mu)\,e_{\mu}] + \\
& + & \Sigma_{\lambda}[d \omega_a^{\lambda} \,e_{\lambda} -
\omega_a^{\lambda}\wedge (\Sigma_c \, \omega_{\lambda}^c\, e_c +
\Sigma_{\mu} \, \omega_{\lambda}^{\mu}\, e_{\mu})],
\end{eqnarray*}
one has \begin{equation} \label{curvaS} \ccR_{(\fE, h)}(e_a) =
\ccR_{(\fS, h)} (e_a) + (B \wedge A) (e_a) + (\D_{(\fS,
h)}A)(e_a),
\end{equation}
\begin{equation} \label{curvaQ} \ccR_{(\fE, h)}(e_{\lambda}) =
\ccR_{(\fQ, h_Q)} (e_{\lambda}) + (A \wedge B) (e_{\lambda}) +
(\D_{(\fQ, h)}B)(e_{\lambda}).
\end{equation}

One may note that the forms $A$ and $B$ are the same as in the
ordinary case, i.e., they do not carry contributions from the
Higgs fields.

\subsection{Approximate Hermitian-Yang-Mills structure and semistability}
For later use, we extend to the Higgs case the relation between
the existence of approximate Her\-mi\-tian-Yang-Mills structure
and semistability. Given an Hermitian vector bundle $(E,h)$,
we introduce a norm on the space of Hermitian endomorphisms $\psi$
of $(E,h)$ by letting
$$ \vert \psi \vert = \max_X \sqrt{\mbox{tr}( \psi^2) }\,.$$

\begin{defin} We say that a Higgs   bundle $\fE=(E, \phi)$
has an \emph{approximate Her\-mi\-tian-Yang-Mills struc\-tu\-re}
if for every positive real number $\xi$ there is an Hermitian
metric $h_{\xi}$ on $E$ such that
\begin{equation}
\vert \mathcal{K}_{(\fE,h)}-c \cdot \Id _E\vert  < \xi \,.
\end{equation}
\end{defin}
The constant $c$ is again given by equation \eqref{c}.

We want to prove the following result (this is proved in the case of vector bundles
in \cite[VI.10.13]{koba}).

\begin{thm} \label{ss} \label{approHYM} A Higgs   bundle $\fE=(E, \phi)$
on a compact K\"ahler manifold
admitting an approximate Her\-mi\-tian-Yang-Mills struc\-tu\-re is
semistable.
\end{thm}

We shall need as usual a vanishing result.

\begin{defin} A section $s$ of a Higgs   bundle $\fE=(E, \phi)$
is said to be \emph{$\phi$-invariant} if it is an eigenvector of
$\phi$, that is, if there is a holomorphic 1-form $\lambda$ on $X$
such that $\phi(s)=\lambda\otimes s$.
\end{defin}

\begin{prop}\label{sezioneparallele}
Assume that a Higgs   bundle $\fE=(E, \phi)$ has an Hermitian
metric $h$ such that the mean curvature $\mathcal{K}_{(\fE, h)}$
of the Hitchin-Simpson connection defines a seminegative definite form.
Then $D_{(E, h)}(s) =0 $ and $\mathcal{K}_{(\fE, h)}(s, s) =0$ for
every $\phi$-invariant section $s$ of $E$.
\end{prop}

\begin{proof}
We start by writing  the relation between the curvatures of the
Chern and Hitchin-Simpson connections for $(\fE, h)$. One has
\begin{equation}\label{curvasimpsochern}
\ccR_{(\fE, h)} = R_{(E, h)} + D'_{(E, h)}(\phi) + D''_{(E,
h)}(\bar\phi) + [\phi, \bar\phi]
\end{equation}
where we have split the Chern connection $D_{(E, h)}=D'_{(E,
h)}+D''_{(E, h)}$ into its (1,0) and (0,1) parts, and $[\phi,
\bar\phi] = \phi\circ\bar\phi + \bar\phi\circ\phi$.

We now compute the Hitchin-Simpson curvature on a $\phi$-invariant section
$s$ of $E$. Writing $\phi(s)=\lambda\otimes s$, one has
$$[\phi, \bar\phi] (s)=0\,,\qquad D'_{(E, h)}(\phi)(s)=\partial \lambda
\otimes s\,, \qquad D''_{(E, h)}(\bar\phi)(s)=\bar\partial\bar
\lambda\otimes s\,.$$ Thus we have
$$
\ccR_{(\fE, h)}(s) =  R_{(E,h)} (s) +
d(\lambda+\bar\lambda)\otimes s\,.$$ Still considering a
$\phi$-invariant section $s$, one has from here and from the
Weitzenb\"ock formula \cite{Besse87} the identity
$$
\partial\bar\partial  h(s, s)) =  h \big(D'_{(E,h)}(s),
D'_{(E,h)}(s) \big) - h \big(\ccR_{(\fE, h)}(s), s \big) +
 h(s, s) \,d(\lambda+\bar\lambda)\,.
$$
Let us set $f=h(s,s)$ and $L(f)=\Lambda(\partial\bar\partial
\,f)$. By applying the operator $\Lambda$ to the previous
equation, we obtain, due to the hypotheses of this Proposition
(and to the fact that the 2-form $d(\lambda+\bar\lambda)$ has no
(1,1) part)
$$L(f) =   \Vert D'_{(E, h)}(s) \Vert ^2 - \mathcal{K}_{(\fE, h)}(s, s)
\geq 0$$
where $\Vert D'_{(E, h)}(s) \Vert ^2 $ is the scalar product
of $D'_{(E, h)}(s) $ with itself using the fibre metric $h$ and the K\"ahler metric on $X$.
By Hopf's maximum principle (see e.g.~\cite{koba}) this implies
$L(f)=0$, which in turn implies $D'_{(E, h)}(s)=0$ and
$\mathcal{K}_{(\fE, h)}(s, s) =0$. Since $s$ is holomorphic, and
the Chern connection is compatible with the holomorphic structure
of $E$, we also have $D_{(E, h)}(s)=0$.
\end{proof}

\begin{corol}\label{sezione}
Let $(\fE, h)$ be an Hermitian Higgs   bundle. If the mean
curvature $\K_{(\fE, h)}$ of the Hitchin-Simpson connection is
seminegative definite everywhere, and negative  definite at some
point of $X$, then $E$ has no nonzero $\phi$-invariant sections.
\end{corol}

\begin{proof}  If $s$ is a nonzero $\phi$-invariant  section of $E$,
then it never vanishes on $X$ since $D_{(E,h)}(s)=0$ by
Proposition \ref{sezioneparallele}. By the same Proposition
$\K_{(\fE, h)}(s,s)=0$, and this contradicts the fact that
$\K_{(\fE, h)}$ is negative at some point.
\end{proof}

\begin{corol} \label{negdeg}
Let $\fE=(E, \phi)$ be a Higgs bundle over $X$ which admits an
approximate Hermitian-Yang-Mills structure. If  $\deg(E) < 0$ then
$E$ has no nonzero  $\phi$-invariant sections.
\end{corol}
\begin{proof}
Since $\fE$ admits an approximate Hermitian-Yang-Mills
structure, for every $\xi > 0$  there exists an Hermitian metric
$h_{\xi}$ on $E$ such that  $$- \xi \cdot h_{\xi} < \K_{(\fE,
h_{\xi})} - c \cdot h_{\xi} < \xi \cdot h_{\xi}$$ with $c<0$. Then
for $\xi$ small enough $\K_{(\fE, h_{\xi})}$ is negative definite,
and the result follows from the previous corollary.
\end{proof}

\noindent\emph{Proof of Theorem \ref{ss}}. Assume that $\fE$
admits an approximate Hermitian-Yang-Mills structure and let $\fF$
be a Higgs subsheaf of $\fE$, with $\rk(F)=p$. Let $\fG$
be the Higgs bundle $(G,\psi)$, where $G= \wedge^pE \otimes
\det(F)^{-1}$, and $\psi$ is the Higgs field naturally induced on
$G$ by the Higgs fields of $\fE$ and $\fF$.  The inclusion
$\fF\hookrightarrow \fE$ induces a morphism $\det (\fF) \to
\wedge^p \fE$, and, tensoring by $\det(\fF)^{-1}$, we obtain a
$\psi$-invariant section of $G$. Now, it is not difficult to check
that the Higgs bundle $\fG$ admits an approximate
Hermitian-Yang-Mills structure, with constant
$$ c_G = \frac{2np\pi}{n!\vol(X)} (\mu(E) -\mu(F))\,.$$
By Corollary \ref{negdeg} we have $c_G\ge 0$, so that $\fE$ is
semistable.

\section{Metric characterization of numerical effectiveness for Higgs
bundles}\label{3}
The usual definitions of numerically effective line and vector bundle,
which are given in terms of embedded curves,
are not appropriate in
the K\"ahler case  since a K\"ahler manifold need not
contain curves. The approaches by Demailly, Peternell and
Schneider \cite{DPS94} and by de Cataldo \cite{dC98} rely on a
definition of numerical effectiveness in terms of fibre metrics.
In particular, de Cataldo considers metrics on the bundle $E$, and
this approach seems to be well suited to an extension to the case
of Higgs bundles, again implementing the idea that the transition
from the ordinary to the Higgs case is obtained by replacing  the
Chern connection with the Hitchin-Simpson connection.

\subsection{Numerically effective Higgs bundles}
De Cataldo's formalism rests on the following terminology. Let
$V$, $W$ be finite-di\-men\-sional complex vector spaces, and let
$\theta_1$, $\theta_2$ be Hermitian forms on $V \otimes W$. Let
$t$ be any positive integer; one writes $\theta_1 \geq_t \theta_2$
if the Hermitian form $\theta_1 - \theta_2$ is semipositive
definite on all tensors in $V \otimes W$ of rank $\rho \leq t$
(where the rank is that of the associated linear map $V^\ast \to
W$). Of course the relevant range for $t$ is $1 \le t \le N =
\min(\dim V,\dim W)$.

If $X$ is a compact K\"ahler manifold of dimension $n$ and $(E,h)$
is a rank $r$ Hermitian vector bundle on $X$, equipped with a
connection $D$, we may associate with the curvature $R$ of $D$ an
Hermitian form $\widetilde{R}$ on $T_X\otimes E$, defined by
\begin{equation} \label{curvatilde} \widetilde{R}(u\otimes s,v\otimes t)
=\frac{i}{2\pi} \langle h({R}^{(1,1)}(s),t) , u\otimes v\rangle
\,.\end{equation} where ${R}^{(1,1)}$ is the $(1,1)$ part of $R$.
According to de Cataldo, the Hermitian bundle $(E, h)$ is $t$-nef
if for every $\xi > 0$ there exists an Hermitian metric $h_{\xi}$
on $E$ such that $\widetilde{R}_ {(E, h_{\xi})} \ge_t -\xi \omega
\otimes h _\xi$.
\begin{defin} A Higgs bundle $\fE=(E, \phi)$ on $X$ is said to be
\begin{enumerate} \item $t$-H-semipositive if there is an Hermitian
metric $h$ on $E$ such that $\widetilde{\ccR}_{(\fE,h)} \geq_t 0$;
\item $t$-H-nef if for every $\xi > 0$ there exists an Hermitian
metric $h_{\xi}$ on $E$ such that $\widetilde{\ccR}_ {(\fE,
h_{\xi})} \ge_t -\xi \omega \otimes h_\xi$; \item $t$-H-nflat if
both $\fE$ and $\fE^*$ are   $t$-H-nef.
\end{enumerate}
\end{defin}

\begin{remark}
\begin{itemize}
    \item [(i)] If $\fE$ is
    $t$-H-semipositive, then $\fE$ is $t$-H-nef. A $t$-H-semipositive
    ($t$-H-nef, respectively) Higgs bundle $\fE$  is
    $t'$-H-semipositive ($t'$-H-nef, respectively) for every $t'$ such
    that $1 \leq t' \leq t$;
    \item [(ii)] since the $(1,1)$ part of the Hitchin-Simpson curvature of a
    Higgs line bundle coincides with the Chern curvature, a Higgs line
    bundle is 1-H-nef if and only if it is 1-nef (as an ordinary
    bundle). More generally, since $\widetilde{\ccR}_{(\fE,h)}
    \geq_t \widetilde{R}_{(E,h)}$, if $E$ is $t$-nef then
    $\fE=(E,\phi)$ is $t$-H-nef for every choice of the Higgs
    field $\phi$;
   \item [(iii)] an Hermitian flat Higgs bundle is $t$-H-nflat for every $t$.
\end{itemize}
\end{remark}

In the next Propositions we establish some basic properties of
$t$-H-nef Higgs bundles on a compact K\"ahler manifold $X$.

\begin{prop} \label{pullback}
Let $f: X \to Y$ be a holomorphic map, where $X$ and $Y$ are
compact K\"ahler manifolds, and let $\fE =(E, \phi)$ be a $t$-H-nef
Higgs vector bundle on $Y$; then $f^{\ast} \fE =(f^{\ast}E,
f^{\ast}\phi)$ is a 1-H-nef Higgs bundle over $X$.
\end{prop}

\begin{proof}
This is proved as in \cite[Proposition 3.2.1(1)]{dC98}.
\end{proof}

\begin{prop} \label{tensor}
Let $\fE =(E, \phi_E)$ and $\fF = (F, \phi_F)$ be Higgs bundles.
If $\fE$ is $t'$-H-nef and $\fF$ is $t''$-H-nef, then $\fE \otimes
\fF =(E \otimes F, \rho)$ is $t$-H-nef, where
\begin{equation}
\begin{array}{cccc}
  \rho: & E \otimes F & \longrightarrow & E \otimes F \otimes
  \Omega^1 \\
  & \rho(e \otimes f) & \mapsto & \phi_{\fE}(e) \otimes f + e \otimes
  \phi_{\fF}(f) \\
\end{array}
\end{equation}
and $t= \min(t',t'')$.
\end{prop}

\begin{proof}
Since $\fE =(E, \phi_E)$ (analog., $\fF = (F, \phi_F)$) is
$t'$-H-nef, for all $\xi > 0$ there exists a metric $h_{(\fE,
\xi/2)}$ over $\fE$ (analog. $h_{(\fF,\xi/2)}$ over $\fF$) such
that the Hitchin-Simpson curvature
$\widetilde{\ccR}_{(\fE,h_{(\fE,\xi/2)})}$ satisfies
$\widetilde{\ccR}_{(\fE,h_{(\fE,\xi/2)})} \geq_{t'} -\frac{\xi}{2}
\omega \otimes h_{(\fE,\xi/2)}$ (analogously,
$\widetilde{\ccR}_{(\fF, h_{(\fF,\xi/2)})} \geq_{t''}
-\frac{\xi}{2} \omega \otimes h_{(\fF, \xi/2)}$). Considering on
$\fE \otimes \fF$ the metrics $h_{\xi} = h_{(\fE,h_{(\fE,\xi/2)})}
\otimes h_{(\fF,h_{(\fF,\xi/2)})}$ we have
$$\widetilde{\ccR}_{(\fE \otimes \fF, h_{\xi})} =
\widetilde{\ccR}_{(\fE, h_{(\fE,{\xi}/{2})})} \otimes h
_{(\fF,\xi/2)} + h _{(\fE,\xi/2)} \otimes \widetilde{\ccR}_{(\fF,
h_{(\fF,{\xi}/{2})})} \geq_t -\xi \omega \otimes h _{\xi}.$$
\end{proof}

Analogously, one proves:

\begin{prop} \label{det}
If $\fE =(E, \phi)$ is a $t$-H-nef Higgs bundle, then for all $p=
2, \ldots, r \, = \rk(E)$ the $p$-th exterior power $\wedge^p \fE
= (\wedge^p E, \wedge^p \phi)$ is a $t$-H-nef Higgs bundle, and
for all $m$, the $m$-th symmetric power $S^m \fE = (S^m E, S^m
\phi)$ is a $t$-H-nef Higgs bundle.
\end{prop}

\begin{lemma} \label{curvatura}
Let $(\fQ, h_{Q})$ be an Hermitian Higgs quotient of $(\fE,h)$.
The respective Hitchin-Simpson curvatures verify the inequality
$\widetilde{\ccR}_{(\fQ,h_{\fQ})} \geq_1
\widetilde{\ccR}_{(\fE,h)\vert \fQ}$.
\end{lemma}

\begin{proof}
The Gauss-Codazzi equations \eqref{curvaQ} imply that the Hitchin-Simpson
curvature of $\fQ$ is given by the restriction of the Hitchin-Simpson
curvature of $\fE$ to $\fQ$ (if we embed $\fQ$ into $\fE$ by
orthogonally splitting the latter) plus the semipositive term $A
\wedge A^{\ast}$, where $A^{\ast}$ is the Hermitian conjugate of
$A$.
\end{proof}

\begin{prop} \label{quot}
A Higgs quotient $\fQ$ of a 1-H-nef Higgs bundle $\fE = (E, \phi)$
is 1-H-nef.
\end{prop}

\begin{proof} Let $\xi > 0$ and $h_{\xi}$ be an Hermitian metric on $\fE$ with
    $\widetilde{\ccR}_ {(\fE, h_{\xi})} \geq_1 -\xi \omega \otimes h _{\xi}$.
    We can endow $\fQ$ with the quotient metric $h_{(\fQ,\xi)}$ and embed
    it into $\fE$ as a $C^\infty$ Higgs subbundle. The claim follows from
    Lemma \ref{curvatura}.
\end{proof}

\begin{corol} \label{kernel}
If $0 \to \fS \to \fE \to \fQ \to 0$ is an exact sequence of Higgs
bundles, with $\fE$ and $\det(\fQ)^{-1}$ 1-H-nef, then $\fS$ is
1-H-nef.
\end{corol}

\begin{proof}
The proof is as in Proposition 1.15(iii) of \cite{DPS94}. Let $r =
\rk(E)$ and $p=\rk(S)$. By taking the $(p-1)$-th exterior power of
the morphism $\fE^{\ast} \to \fS^{\ast}$ obtained from the exact
sequence in the statement, and using the isomorphism $\fS \simeq
\wedge^{p-1} \fS^{\ast} \otimes \det(\fE)$, we get a surjection
$\wedge^{p-1} \fE^{\ast} \to \fS \otimes \det(\fS)^{-1}$.
Tensoring by $\det(\fS) \simeq \det(\fE) \otimes \det(\fQ)^{-1}$
we eventually obtain a surjection $\wedge^{r-p+1} \fE \otimes
\det(\fQ)^{-1} \to \fS$. Propositions \ref{tensor} and \ref{quot}
now imply the claim.
\end{proof}

\begin{prop} \label{extension}
An extension of 1-H-nef Higgs bundles is 1-H-nef.
\end{prop}

\begin{proof}
Let us consider an extension of Higgs bundles
$$0 \to \fF \to \fE \to \fG \to 0$$ where $\fF$ and $\fG$ are 1-H-nef.
Then for every $\xi > 0$ the latter bundles carry Hermitian
metrics $h_{(\fF,\xi)}$ and $h_{(\fG,\xi)}$ such that
$$\widetilde{\ccR}_{(\fF,h_{(\fF,\xi/3)})} \ge_1 - \tfrac{\xi}{3}
\omega \otimes h_{(\fF,\xi/3)}, \qquad
\widetilde{\ccR}_{(\fG,h_{(\fG,\xi/3)})} \ge_1 - \tfrac{\xi}{3}
\omega \otimes h_{(\fG,\xi/3)}.$$ Fixing a $C^\infty$  isomorphism $\fE
\simeq \fF \oplus \fG$, these metrics induce an Hermitian metric
$h_{\xi}$ on $\fE$. A simple calculation, which involves the
second fundamental form of $\fE$, shows that
$\widetilde{\ccR}_{(\fE,h_{\xi})} \ge_1 - \xi \omega \otimes
h_\xi$, so that $\fE$ is 1-H-nef (details in the ordinary case are
given in \cite{dC98}).
\end{proof}

\subsection{Numerical effectiveness and semistability}
We study here the relations between the metric characterization of
the numerical effectiveness of Higgs bundles and their
semi\-stability. The following result will be a useful tool.
\begin{prop} \label{sectdual}
Let $\fE=(E, \phi)$ be a 1-H-nef Higgs bundle, and $\fE^{\ast}
=(E^{\ast}, \phi^{\ast})$ the dual Higgs bundle. If $s$ is a
$\phi^{\ast}$-invariant section of $E^{\ast}$, then $s$ has no
zeroes.
\end{prop}

\begin{proof}
We modify the proof of Proposition 1.16 in \cite{DPS94}. For a given $\xi >
0$, let $h_{\xi}$ be the metric on $E$ such that
$\widetilde{\ccR}_ {(\fE, h_{\xi})} \ge_1 -\xi \omega \otimes
h_\xi$. Let us define the closed $(1,1)$ current $$T_{\xi} =
\frac{i}{2\pi} \partial \bar \partial \log h^{\ast}_{\xi}(s,s);$$
a simple computation shows that it satisfies the inequality
$$T_{\xi} \geq -\frac{\widetilde{R}_ {(E^{\ast},
h^{\ast}_{\xi})}(s,s)}{h^{\ast}_{\xi}(s,s)},$$ where
$\widetilde{R}_ {(E^{\ast}, h^{\ast}_{\xi})}(s,s)$ is regarded as
a 2-form on $X$. Now, if $s$ is $\phi^{\ast}$-invariant, then
$[\phi^{\ast}, \overline{\phi}^{\ast}](s)=0$, so that
$\widetilde{\ccR}_ {(\fE^{\ast}, h^{\ast}_{\xi})}(s,s) =
\widetilde{R}_ {(E^{\ast}, h^{\ast}_{\xi})}(s,s)$. On the other
hand, since $\fE$ is 1-H-nef, we have $-\widetilde{\ccR}_
{(\fE^{\ast}, h^{\ast}_{\xi})}(s,s) \geq - \xi \, h_\xi(s,s) \,
\omega.$ Thus, $T_{\xi} \geq - \xi \omega$.

Since $\partial \bar \partial \omega^{n-1}=0$, we have
$$\int_X (T_{\xi}+ \xi \omega) \wedge \omega^{n-1}= \xi \int_X
\omega^n.$$ For $\xi$ ranging in the interval $(0,1]$ the masses
of the currents $T_\xi + \xi \omega$ are uniformly bounded from
above, so that the sequence $\{T_\xi + \xi \omega\}$ contains a
subsequence which, by weak compactness, converges weakly to zero.
(For details on this technique see e.g.~\cite{D01}). Therefore,
the Lelong number of $T_{\xi}$ at each point $x \in X$ (which
coincides with the vanishing order of $s$ at that point) is zero
\cite{S74}, which implies that $s$ has no zeroes.
\end{proof}

\begin{thm} \label{nfss}
A 1-H-nflat Higgs bundle $\fE = (E, \phi)$ is semistable.
\end{thm}

\begin{proof}
Since $\fE$ is 1-H-nef for every $\xi > 0$ it carries an Hermitian
metric $h_{\xi}$ such that $\widetilde{\ccR}_{(\fE, h_{\xi})}
\geq_1 -\xi \omega \otimes h_\xi$. As the mean curvature
$\K_{(\fE, h_{\xi})}$ may be written in the form
$$\K_{(\fE, h_{\xi})}(s,s) = -2\pi \sum_{i=1}^n
\widetilde{\ccR}_{(\fE, h_{\xi})} (e_i \otimes s, e_i \otimes
s),$$ where the $e_i$'s are a unitary frame field on $X$, one has
$$\K_{(\fE, h_{\xi})}(s,s) \leq 2\pi\,n \, \xi \, h_{\xi}(s,s).$$

On the other hand, since $\det(\fE)^{-1}$ is 1-H-nef, the Higgs
bundle $\fE^{\ast} \simeq \wedge^{r-1}\fE \otimes \det(\fE)^{-1}$
is 1-H-nef with the dual metric $h_{\xi}^*$, so that $\K_{(\fE^*,
h_{\xi}^*)} = -\K_{(\fE, h_{\xi})}^{t}$, and
\begin{equation*}
\K_{(\fE, h_{\xi})}(s,s) \geq -2\pi\,n \, \xi \, h_{\xi}(s,s).
\end{equation*}
As $c_1(E)=0$ because $\det(\fE)$ is 1-nflat \cite[Corollary
1.5]{DPS94}, after rescaling $\xi$ these equations imply $\vert
\K_{(\fE, h_{\xi})} \vert \leq \xi $, so that $\fE$ is semistable
by Theorem \ref{approHYM}.
\end{proof}

\subsection{A characterization of 1-H-nflat Higgs bundles}
The next lemmas will be used in the proof of Theorem
\ref{salamancasorbona}, which is one of the main results in this
paper.
\begin{lemma}\label{c10}
A 1-H-nef Higgs bundle $\fE = (E, \phi)$ such that $c_1(E)=0$ is
1-H-nflat.
\end{lemma}

\begin{proof}
This follows again from the fact that $\fE^* \simeq \wedge^{r-1}
\fE \otimes \det(\fE)^{-1}$ is an isomorphism of Higgs bundles.
\end{proof}

\begin{lemma} \label{c10line}
A 1-H-nef Higgs line bundle $\fL$ of zero degree is Hermitian
flat.
\end{lemma}

\begin{proof}
This is already contained in \cite[Cor.~1.19]{DPS94}, but for the
reader's convenience we give here a proof. For every $\xi
> 0$ one has on $\fL$ an Hermitian metric $k_{\xi}$ satisfying the
inequality
\begin{equation*} 0 \leq \int_X
(\tfrac{i}{2\pi}\ccR_{(\fL,k_{\xi})} + \xi \omega) \cdot
\omega^{n-1} = \deg(L) + \xi \int_X \omega^n.
\end{equation*}
By the same argument as in the proof of Proposition
\ref{sectdual}, if $\deg(L)=0$ by taking the limit $\xi \to 0$ one
shows that $c_1(L)=0$, so that $\fL$ is Hermitian flat.
\end{proof}

\begin{lemma}\label{new} If the Higgs bundle $\fE$ is 1-Hnflat,
and $\{h_\xi\}$ is a family of metrics which makes $\fE$ 1-H-nef,
then the family of dual metrics $\{h^\ast_\xi\}$ makes
$\fE^\ast$ 1-H-nef.
\end{lemma}
\begin{proof} The determinant line bundle $\det(E)$ is 1-nef
with respect to the family of determinant metrics $\{\det h_\xi\}$.
The dual line bundle $\det^{-1}(E)$ is 1-nef as well, and
it is such with respect to a family $\{a(\xi)\,\det^{-1} h_\xi\}$,
where the homothety factor $a(\xi)$ only depends on $\xi$ \cite[Cor.~1.5]{DPS94}.
From the isomorphism $\fE^\ast \simeq \wedge^{r-1}\fE\otimes  \det^{-1}(\fE)$
(where $r=\rk(E)$) we see that $\fE^\ast$ is made 1-H-nef
by the family of metrics $\{h'_\xi = a(\xi)\,h^\ast_\xi\}$, so that
for every $\xi>0$ the condition
$\widetilde{\ccR}_{(\fE^\ast,h'_\xi)} \ge_1 - \xi\,\omega\otimes h'_\xi$ holds.
But this implies $\widetilde{\ccR}_{(\fE^\ast,h^\ast_\xi)} \ge_1 - \xi\,\omega\otimes h^\ast_\xi$.
\end{proof}

\begin{lemma} \label{c20}
A stable 1-H-nflat Higgs bundle $\fE$ is Hermitian flat
(cf.~Definition \ref{simpsonconnection}).
\end{lemma}

\begin{proof}
As before, let us denote by $\Vert \ccR_{(\fE,h_\xi)} \Vert^2$ the scalar product
of the Hitchin-Simpson curvature with itself
obtained by using the Hermitian metric of the bundle $E$
and the K\"ahler form on $X$ (thus,  $\Vert
\ccR_{(\fE,h_\xi)} \Vert$ is a function on $X$). Note that
in terms of a local orthonormal frame $\{e_\alpha\}$ on
$X$ and a local orthonormal basis of sections $\{s_a\}$ of $E$
we may write
$$\Vert
\ccR_{(\fE,h_\xi)} \Vert^2 = 4\pi^2 \sum_{\alpha,a}
\left( \widetilde{\ccR}_{(\fE,h_\xi)}(e_\alpha\otimes s_a, e_\alpha\otimes s_a) \right)^2\,.$$
Since $\fE$ is 1-H-nef, for every $\xi>0$ there
is an Hermitian metric $h_\xi$ on $E$ such that $ \widetilde{\ccR}_{(\fE,h_\xi)} \ge_1 - \xi \,\omega\otimes h_\xi$. Taking Lemma \ref{new} into account, for every $\xi$ we have the
inequalities
$$ \xi \ge \widetilde{\ccR}_{(\fE,h_\xi)}(e_\alpha\otimes s_a, e_\alpha\otimes s_a) \ge - \xi\,.$$
So we have $\Vert \ccR_{(\fE,h_\xi)} \Vert \le a_1\,\xi$ for some constant $a_1$.
In the same way we have  $\Vert \mathcal K_{(\fE,h_\xi)} \Vert \le a_2\,\xi$ for some constant $a_2$.

Assume that where $n=\dim X>1$.
Since $c_1(E)=0$, we have the representation formula \cite[Chap. IV.4]{koba}
\begin{equation*}
\int_X c_2(E) \cdot \omega^{n-2} = \frac{1}{8\pi^2\,n(n-1)} \int_X
(\Vert \ccR_{(\fE,h_\xi)}
\Vert^2 - \Vert \K_{(\fE,h_\xi)} \Vert^2) \, \omega^n
\end{equation*}

The previous inequalities imply $\int_X c_2(E) \cdot
\omega^{n-2} = 0$. For every value of $n$,  may apply Theorem 1 and
Proposition 3.4 in \cite{S88} to show that $\fE$ admits an
Hermitian metric whose corresponding Hitchin-Simpson curvature vanishes.
\end{proof}

\begin{thm} \label{salamancasorbona} A Higgs bundle $\fE$
is   1-H-nflat if and only if it has a filtration in Higgs subbundles whose quotients
are Hermitian flat Higgs bundles. As a consequence,  all Chern
classes of a 1-H-nflat Higgs bundle vanish.
\end{thm}

\begin{proof} Assume that $\fE$ has such a filtration. Then
any quotient of the filtration is 1-H-nflat, and the claim follows
from Proposition \ref{extension}.

To prove the converse,
let $\fF$ be a Higgs subsheaf of $\fE$ of rank $p$. We have an
exact sequence of Higgs sheaves $$0 \to \det(\fF) \to \wedge^p \fE
\to \fG \to 0,$$ where $\fG$ is not necessarily locally-free. Since $\det(\fE)$ is
1-H-nflat we have $c_1(E)=0$. By
Theorem \ref{nfss} $\wedge^p \fE$ is semistable, so that $\deg(F)
\leq 0$. Let $h_{\xi}$ be a family of Hermitian metrics which
makes $\fE$ a 1-H-nef Higgs bundle, and let $h^p_{\xi}$ be the
induced metrics on $\wedge^p\fE$. After rescaling the dual metrics
$(h^p_{\xi})^{\ast}$ we obtain a family of metrics which makes
$\wedge^p \fE^{\ast}$ a 1-H-nef Higgs bundle (cf.~Lemma \ref{new}). Let $U$ be the open
dense subset of $X$ where $\fG$ is locally free; then the metrics
$(h^p_{\xi})^{\ast}$ induce on $\det(\fF)^{-1}_{\vert U}$ metrics
making it 1-H-nef. These metrics extend to the whole of $X$, since
they are homothetic by a constant factor to the duals of the
metrics induced on $\det(\fF)$ by the metrics on $\wedge^p \fE$.
Thus, $\det(\fF)^{-1}$ is 1-H-nef.
If $\deg(F)=0$ by Lemma \ref{c10line} $\det(\fF)$ is Hermitian
flat, so that $\wedge^p \fE \otimes \det(\fF)^{-1}$ is 1-H-nflat.
Then by Proposition \ref{sectdual} the morphism of Higgs bundles
$\det(\fF) \to \wedge^p \fE$ has no zeroes, so that $\fG$ is
locally-free.

In view of Lemma \ref{c20} we may assume that $\fE$ is not stable.
Let us then identify $\fF$ with a destabilizing Higgs subsheaf of
minimal rank and zero degree. We need $\fF$ to be reflexive;
we may achieve this by replacing $\fF$ with
its double dual $\fF^{\ast \ast}$. By Lemma 1.20 in \cite{DPS94},
$\fF$ is locally-free and a Higgs subbundle of $\fE$. Now,
$\fF^{\ast}$ is 1-H-nef because it is a Higgs quotient of
$\fE^{\ast}$, while $\fF$ is 1-H-nef by Corollary \ref{kernel}, so
that $\fF$ is 1-H-nflat. Since $\fF$ is stable by construction, by
Lemma \ref{c20} it is Hermitian flat. The existence of the
filtration follows by induction on the rank of $\fE$
since the quotient $\fE/\fF$ is locally-free and 1-H-nflat, hence
we may apply to it the inductive hypothesis.
\end{proof}

\subsection{Projective curves} We conclude this section by proving
two results that hold when $X$ is a smooth projective curve. The
first Proposition generalizes results given in
\cite{G79,BHR05,BG05}.
\begin{prop} \label{gieseker} If a Higgs bundle $\fE$ on $X$ is
semistable and $\deg(E) \ge 0 $, then $\fE$  is  1-H-nef.
\end{prop}
\begin{proof} If $\fE$ is stable it admits an Hermitian-Yang-Mills
metric $h$, so that $ \widetilde{\ccR}_{(\fE,h)} = c\, h$
with $c \ge 0$ (note that we essentialy identify $\widetilde{\ccR}_{(\fE,h_{\xi})} $
with the mean curvature since we are on a curve). Then $\fE$  is  1-H-nef.

If $\fE$ is properly semistable, we may filter it in such a way that the quotients
of the filtration are  stable
Higgs bundles of nonnegative degree. By the previous argument, every
quotient is 1-H-nef. One then concludes by Proposition \ref{extension}.
\end{proof}

In \cite{H66} Hartshorne gives a characterization of ample vector bundles
of rank 2 on a smooth projective curve. Our next result generalizes
this to 1-H-nef Higgs bundles of any rank on a smooth projective curve.
A similar statement might be easily obtained for the ample case. This
also partly generalizes Theorem 3.3.1 in \cite{dC98}.

\begin{prop} \label{quotient}
Let $\fE$ be a Higgs bundle of nonnegative degree over a smooth
projective curve $X$, whose locally-free Higgs quotients are all
1-H-nef. Then $\fE$ is 1-H-nef.
\end{prop}

\begin{proof} One may consider two cases:

1.- If $\fE$ is semistable then by Proposition \ref{gieseker}
it is 1-H-nef.

2.- If $\fE$ is not semistable, then it has a semistable Higgs
subbundle $\fK$ with $\mu (K) > \mu(E)$, so that one has a exact
sequence of Higgs bundles $0 \rightarrow \fK \rightarrow \fE
\rightarrow \fQ \rightarrow 0$. Since $\deg(E) \geq 0$, then
$\deg(K) > 0$, and again by Proposition \ref{gieseker} we have
that $\fK$ is 1-H-nef. Thus $\fE$ is an extension of 1-H-nef Higgs
bundles, so that it is 1-H-nef by Proposition \ref{extension}.
\end{proof}

\begin{remark} It seems useful to   state in a close way the relations between
the conditions of 1-H-nefness and semistability in the case of Higgs bundles
on curves. So, let $X$ be a complex smooth projective curve, and
$\fE=(E,\phi)$ a Higgs bundle on $X$. Let $d=\deg(E)$.
\begin{enumerate} \item If $\fE$ is semistable and $d>0$, then
$\fE$ is 1-H-nef.
\item If $\fE$ is semistable and $d<0$, then $\fE^\ast$ is 1-H-nef.
\item if $\fE$ is semistable and $d=0$, then $\fE$ is 1-H-nflat.
\end{enumerate}
On the other hand, if $\fE$ is 1-H-nef and $d\ne 0$, then it need
not be semistable (an example is provided by a direct sum of
1-H-nef Higgs bundles of different slopes). However, if $\fE$ is
1-H-nef and $d=0$ (i.e., it is 1-H-nflat) we know it is semistable
(Theorem \ref{nfss}; this also follows from Corollary 3.6 of
\cite{BG05} since, as we shall see in the next Section, on a curve
the notions of H-nefness and 1-H-nefness coincide).
\end{remark}

\section{The projective case} \label{4}
In our previous paper \cite{BG05} we gave a definition of
numerical flatness for Higgs bundles on smooth projective
varieties. In this section we want to compare it with the
definition we have given here in the case of K\"ahler manifolds.
We start by briefly recapping the situation in the projective
case.

\subsection{Grassmannians of Higgs quotients}
Our definition of H-nefness for Higgs bundles on projective varieties
requires to consider Higgs bundles on  singular  spaces
(the Higgs Grassmannians whose definition we are going to recall
in this Section, see also \cite{BHR05,BG05}). As we shall see, Higgs Grassmannians
are locally defined as the zero locus of a set of holomorphic functions
on the usual Grassmannian varieties, and are therefore schemes.
For such spaces there is well-behaved theory of the de Rham complex
\cite{EGAIV}, which is all one needs to define Higgs bundles. Let us in particular
notice that Higgs bundles on schemes are well-behaved with respect to
restrictions to closed subschemes.
Indeed, if $\fE=(E,\phi)$ is a Higgs
bundle on a scheme $X$, and $Y \hookrightarrow X$ is a closed immersion,
due to the isomorphism
\begin{eqnarray*}(E\otimes_{\cO_X}\Omega_X)_{\vert Y}
&= & (E\otimes_{\cO_X}\Omega_X) \otimes_{\cO_X} \cO_Y  \\
&\simeq & (E \otimes_{\cO_X} \cO_Y ) \otimes_{\cO_Y}  (\Omega_X\otimes_{\cO_X} \cO_Y )
\simeq E _{\vert Y} \otimes_{\cO_Y} {\Omega_X}_{\vert Y}\,, \end{eqnarray*}
by composing the restriction $\phi_{\vert Y} \colon E _{\vert Y} \to
E _{\vert Y} \otimes_{\cO_Y}
{\Omega_X}_{\vert Y}$ with the projection $p\colon{\Omega_X}_{\vert Y}\to  \Omega_Y$
we obtain a Higgs bundle $\fE _{\vert Y}=(E_{\vert Y},(1\times p)\circ \phi_{\vert Y})$
on $Y$.

Moreover,
the Chern classes we shall be using in that case
may be regarded as those in the Fulton-MacPherson theory, cf.~\cite{Ful}.

Thus, let $X$ be a scheme over the complex numbers. Given a rank $r$
vector bundle $E$ on $X$, for every positive integer $s$ less than
$r$, let $\grass_s(E)$ be the Grassmann bundle of $s$-planes in
$E$, with projection $p_s : \grass_s(E) \to X$. This is a
parametrization of the rank $s$ locally-free quotients of $E$.
There is a universal exact sequence
\begin{equation}\label{univ}
0 \to S_{r-s,E} \xrightarrow{\psi} p_s^*(E) \xrightarrow{\eta}
Q_{s,E} \to 0
\end{equation}
of vector bundles on $\grass_s(E)$, with $S_{r-s,E}$ the universal
rank $r-s$ subbundle and $Q_{s,E}$ the universal rank $s$ quotient
bundle \cite{Ful}.

Given a Higgs bundle $\fE $, we may construct closed subschemes
$\hgrass_s(\fE)\subset \grass_s(E)$ pa\-ram\-e\-tr\-iz\-ing rank
$s$ locally-free Higgs quotients. With reference to the exact
sequence eq.~\eqref{univ}, we define $\hgrass_s(\fE)$ as the
closed subscheme
of $\grass_s(E)$ where the composed morphism
$$(\eta\otimes1)\circ p_s^\ast(\phi) \circ \psi\colon S_{r-s,E}\to
Q_{s,E}\otimes
 p_s^\ast(\Omega_X)$$ vanishes. We denote by $\rho_s$ the projection $\hgrass_s(\fE)\to
X$. The restriction of \eqref{univ} to the scheme $\hgrass_s(\fE)$
provides the exact sequence (of vector bundles, for the moment)
\begin{equation} \label{univg}
0 \to S_{r-s,\fE}\to \rho_s^\ast (\fE) \to Q_{s,\fE}\to 0 \,.
\end{equation}
The Higgs morphism $\phi$ of $\fE$ induces by pullback
a Higgs morphism $\Phi\colon  \rho_s^\ast (\fE) \to
 \rho_s^\ast (\fE) \otimes \Omega_{\hgrass_s(\fE)}$
 (note in particular that $\Phi\wedge\Phi=0$).
On the account of the condition
 $(\eta\otimes1)\circ p_s^\ast(\phi) \circ \psi=0$ which holds true
 on $\hgrass_s(\fE)$, the morphism $\Phi$ sends $ S_{r-s,\fE}$ to $ S_{r-s,\fE}\otimes \Omega_{\hgrass_s(\fE)}$.  Therefore, $S_{r-s,\fE}$ is a Higgs subbundle
 of $\rho_s^\ast (\fE)$, and the quotient $ Q_{s,\fE}$ inherits a structure of Higgs bundle
 as well. The sequence \eqref{univg} is therefore an exact sequence
 of Higgs bundles.

 By its very construction, the scheme $\hgrass_s(\fE)$ and the quotient
bundle $ Q_{s,\fE}$ enjoy a universal property: for every morphism $f\colon Y\to X$ and every  rank
$s$ Higgs quotient $\fF$ of $f^\ast( \fE)$ there is a morphism
$g\colon Y\to \hgrass_s(\fE)$ such that $f=\rho_s\circ g$ and
$\fF\simeq g^\ast (Q_{s,\fE})$ as Higgs bundles. Therefore, the scheme
 $\hgrass_s(\fE)$ will be called  the \emph{Grassmannian
of locally-free rank $s$ Higgs quotients} of $\fE$, and
$Q_{s,\fE}$ will be called the \emph{rank $s$ universal Higgs  quotient vector
bundle}.

 It is useful
to introduce the following numerical classes
\begin{equation} \label{lambda}
\lambda_{s,\fE}= \left[ c_1(\OPQ{s,\fE}(1) -\frac1r
\pi^\ast_s ( c_1(E))\right]  \in N^1(\PP Q_{s,\fE})
\end{equation}
\begin{equation}\label{theta}
\theta_{s, \fE}=\left[ c_1(Q_{s,\fE})-\frac sr\, \rho_s^\ast (c_1(E)) \right]
\in N^1(\hgrass_s(\fE)), \end{equation} where, for every
projective scheme $Z$, we denote by $N^1(Z)$ the vector space of
$\R$-divisors modulo numerical equivalence: $$N^1(Z) =
\frac{\mbox{Pic}(Z)}{num.\,eq.} \otimes \R.$$

\subsection{Comparison between the projective and K\"ahlerian cases}
We recall from \cite{BG05} the notion of H-nef Higgs bundle.

\begin{defin}  \label{moddef} A Higgs bundle $\fE$ of rank one is said
to be Higgs-numerically effective (for short, H-nef) if it is
numerically effective in the usual sense. If $\rk \fE \geq 2$ we
require that:
\begin{enumerate} \item all bundles $Q_{s,\fE}$ are Higgs-nef;
\item the line bundle $\det(E)$ is nef.
\end{enumerate}
If both $\fE$ and $\fE^\ast$ are Higgs-numerically effective,
$\fE$ is said to be Higgs-numerically flat (H-nflat).
\end{defin}

If the Higgs field is zero (i.e., $\fE$ is an ordinary vector bundle)
this definition reduces to the usual one.

\begin{remark} \label{reduce}
Due to our iterative definition of H-nefness, a Higgs bundle $\fE$ is H-nef if and only if a finite number of line bundles $L_i$ (each defined on a projective scheme $Y_i$ for
which a surjective morphism $Y_i \to X$ exists) are nef. For instance, if
$\fE$ is a rank 3 Higgs bundle of $X$, one is requiring the usual nefness of the following  line bundles:
\begin{itemize} \item  $\det(\fE)$ on $X$ \item $Q_{1,\fE}$ on $\hgrass_1(\fE)$ \item
$\det (Q_{2,\fE})$ on $\hgrass_2(\fE)$  \item
$Q_{1,Q_{2,\fE}}$ on $\hgrass_1(Q_{2,\fE})$.
\end{itemize}
\end{remark}

\begin{prop}
A 1-H-nef Higgs bundle $\fE = (E, \phi)$ is H-nef.
\end{prop}

\begin{proof}
We proceed by induction on the rank $r$ of  $\fE$. If $r = 1$
there is nothing to prove. If $r >1$, for every $s = 1, \dots,
r-1$ let us consider the universal sequence \eqref{univg} on the
Higgs Grassmannian $\hgrass_s(\fE)$.  Since the Higgs Grassmannian
is in general singular, we consider a resolution of singularities
$\beta_s: B_s(\fE) \to \hgrass_s(\fE)$, and pullback the universal
sequence to $B_s(\fE)$: $$0 \to \beta^{\ast}_sS_{r-s,\fE} \to
\gamma^{\ast}_s\fE \to \beta^{\ast}_sQ_{s,\fE} \to 0,$$ where
$\gamma_s = \rho_s \circ \beta_s $. Since $\fE$ is 1-H-nef, the
pullback $\gamma^{\ast}_s(\fE)$ is 1-H-nef as well, and its Higgs
quotient $\beta^{\ast}_s Q_{s,\fE}$ is 1-H-nef, hence H-nef by the
inductive hypothesis.

We need to show that $Q_{s,\fE}$ is H-nef; in view of Remark
\ref{reduce}, by base change this reduces to proving the following
fact: if $f_i:Z_i \to Y_i$ are surjective morphisms of projective
schemes, and $L_i$ are line bundles on $Y_i$ such that the
pullbacks $f_i^{\ast}L_i$ are nef, then the line bundles $L_i$ are
nef. This follows from \cite[Prop.~2.3]{Fu83}.
\end{proof}

\begin{prop} \label{equiv}
A Higgs bundle $\fE = (E, \phi)$ over a smooth
projective curve $X$ is 1-H-nef if and only if it is H-nef.
\end{prop}

\begin{proof}
We have just proved the necessary condition. We prove the
sufficiency again by induction on the rank $r$ of $\fE$. If $r= 1$
there is nothing to prove. If $r > 1$, note that since $\fE$ is
H-nef, then $\deg (E) \ge 0$, and all its quotients $\fQ$ are
H-nef. By the inductive hypothesis, all $\fQ$ are 1-H-nef; one
concludes by Proposition \ref{quotient}.
\end{proof}

This strongly simplifies the proof of Theorem 3.3.1 of
\cite{dC98}, which gives the same result in the case of ordinary
bundles.

We may use these results to  prove some properties of H-nef Higgs
bundles in addition to those given in \cite{BG05}.

\begin{lemma} \label{curves}
 A Higgs bundle  $\fE = (E, \phi)$ over a smooth
projective variety $X$ is H-nef if and only if $\fE_{\vert C} =
(E_{\vert C}, \phi_{\vert C})$ is H-nef for all irreducible curves $C$ in $X$.
\end{lemma}

\begin{proof}
By Remark \ref{reduce} the Higgs bundle $\fE$ is H-nef if and only
if a finite number of line bundles $L_i$ (each defined on a
projective scheme $Y_i$ for which a surjective morphism $Y_i \to
X$ exists) are nef. The claim then follows.
\end{proof}

\begin{remark} We should note that the fact that a Higgs bundle
restricts to a semistable Higgs bundle on any embedded curve is not enough to
ensure that it is H-nef. Consider for instance a  a surface with Picard number 1, and let
$\fE=(E,\phi)$  be a Higgs bundle of negative degree which is semistable after
restriction to every curve (i.e., it is semistable and $\Delta(E)=0$).
The Higgs bundle $\fE$ cannot be H-nef since it has negative degree.
\end{remark}

\begin{prop}
Let $0 \to \fF \to \fE \to \fG \to 0$ be an exact sequence of Higgs bundles
over a smooth projective variety $X$. If $\fF$ and $\fG$ are
H-nef then $\fE$ is H-nef.
\end{prop}

\begin{proof}
In view of Lemma \ref{curves} we may assume that $X$ is a curve.
The result then follows from Propositions \ref{extension} and
\ref{equiv}.
\end{proof}

In the same way, by using Lemma \ref{curves} one can prove that
the tensor, exterior and symmetric products of H-nef Higgs bundles
are H-nef, thus completing the results given in \cite{BG05}.
Moreover we have:

\begin{prop} \label{sym}
Let $\fE$ be a Higgs bundle. If $S^m(\fE)$ is H-nef for some $m$,
then $\fE$ is H-nef.
\end{prop}

\begin{proof}
Since a rank $s$ Higgs quotient of $\fE$ yields a Higgs quotient
of $S^m(\fE)$ of rank $$N_{(m,s)}={m+s-1 \choose s-1},$$ one has a
morphism $g: \hgrass_s(\fE) \to \hgrass_{N_{(m,s)}}(S^m(\fE))$
such that $g^{\ast}(Q_{N_{(m,s)},S^m(\fE)}) \simeq
S^m(Q_{s,\fE})$. Since $S^m(\fE)$ is H-nef, the symmetric product
$S^m(Q_{s,\fE})$ is H-nef. The claim follows by induction on the
rank of $\fE$.
\end{proof}

\subsection{Semistability criteria}
The notion of H-nefness may be used to provide a characterization
of a special class of Higgs bundles, namely, the semistable Higgs
bundes for which the cohomology class
$\Delta(E)=c_2(E)-\frac{r-1}{2r}c_1(E)^2$ vanishes. Recently
several similar criteria have appeared in the literature, which
all generalize Mi\-yaoka's criterion for the semistability of
vector bundles on projective curves \cite{Mi}. In \cite{BHR05} a
criterion was given for characterizing semistable Higgs bundles
with vanishing $\Delta$ class on complex projective manifolds of
any dimension in terms of the nefness of a certain set of
divisorial classes (it is interesting to note that Higgs bundles
of this type have the property of being semistable after
restriction to any curve in the base manifold). Analogous criteria
have been formulated for principal bundles on complex projective
manifolds \cite{BB04}, and more generally for principal bundles on
compact K\"ahler manifolds \cite{BSch05}.

We discuss now three equivalent conditions of this kind. One of
these is stated in terms of the Higgs bundles $T_{s,\fE} =
S^{\ast}_{r-s,\fE} \otimes Q_{s,\fE}$ on the Higgs Grassmannians
$\hgrass_s(\fE)$. For an ordinary vector bundle $E$, the bundle
$T_{s,E} $ is the vertical tangent bundle to $p_s:\grass_{s}(E)
\to X$.

\begin{thm}\label{critproj}
Let $\fE=(E, \phi)$ be a rank $r$ Higgs bundle on a complex
projective manifold $X$. The following three conditions are
equivalent:
\begin{enumerate}
    \item  The Higgs bundle $\fF=S^r(\fE)\otimes  (\det \fE)^{-1}$ is
    H-nflat;
    \item  $\fE$ is semistable and
    $\Delta(E)=c_2(E)-\frac{r-1}{2r}c_1(E)^2=0$;
    \item the Higgs bundles $T_{s,\fE}$ are all H-nef.
\end{enumerate}
\end{thm}

\begin{proof}
We first prove that (i) implies (ii). Since $\det (\fF)$ is
trivial, the dual Higgs bundle $\fF^\ast$ is H-nef as well, i.e.,
$\fF$ is H-nflat, hence semistable by Theorem 3.1 of \cite{BG05}.
Then, the Higgs bundle $\fF\otimes \fF^\ast\simeq S^r(\fE)\otimes
S^r(\fE^\ast)$ is semistable. This implies that $\fE$ is
semistable.

One also has that $S^r(\fE)\otimes S^r(\fE^\ast)$ is H-nflat so
that its Chern classes vanish by Corollary 3.2 of \cite{BG05}. But
since
$$c_2(S^r(E)\otimes S^r(E^\ast)) = 4r (\mbox{rk} \, S^r(E))^2 \Delta(E) $$
we conclude.

(ii) implies (i): we have that $\fF$ is semistable and
\begin{equation*}
c_1(F)=0, \qquad c_2(F)= 2r (\mbox{rk} \, S^r(E))^2 \Delta(E)=0.
\end{equation*}
By Theorem 2 of \cite{S92}, $\fF$ has a filtration whose quotients
are stable Higgs bundles with vanishing Chern classes. Proceeding
as in the proof of Lemma \ref{c20}, these quotients are shown to
be   Hermitian flat, hence they are H-nflat. Then $\fF$ is
H-nflat as well.

We prove now that (i) implies (iii). If $\fF$ is H-nef, then the
$\Q$-Higgs bundle $\fE \otimes (\det(\fE))^{-1/r}$ is H-nef by Proposition
\ref{sym}, so
that the Higgs bundle $\fE^{\ast} \otimes (\det(\fE))^{1/r}$ is H-nef
(since $c_1(\fE \otimes (\det(\fE))^{-1/r})=0$), and
$Q_{s,\fE} \otimes \rho^{\ast}_s(\det(\fE))^{-1/r}$ is H-nef as
well, since it is a universal quotient of $\fE \otimes (\det(\fE))^{-1/r}$.
From the exact sequence $$0 \to Q_{s, \fE} \otimes Q_{s,
\fE}^{\ast} \to \rho^{\ast}_s(\fE^\ast) \otimes Q_{s, \fE} \to T_{s,\fE}
\to 0,$$ one obtains the claim.

Finally, we prove that (iii) implies (ii). Note that the class
$\theta_{s,\fE}$ defined in equation \eqref{theta} equals
$[c_1(T_{s,\fE})]$, so that if $T_{s,\fE}$ is H-nef, the class
$\theta_{s,\fE}$ is nef. This holds true for every $s=1,\dots,
r-1$. It was proved in \cite{BHR05} that this is equivalent to
condition (ii) in the statement.
\end{proof}

\begin{remark} In \cite{BHR05} it was proved that condition (ii) is fulfilled if
and only if all classes $\lambda_{s,\fE}$ are nef. It may be
interesting to check directly that the latter condition is
equivalent to condition (i) in our Theorem \ref{critproj}.

If $\fF = S^r \fE \otimes \det(\fE)^{-1}$ is H-nef, since
$c_1(F)=0$ it is also H-nflat, so that all its Higgs quotients are
nef in the usual sense (see \cite[Cor.~3.6]{BG05}). Since the
Higgs bundle $S^r Q_{s,\fE} \otimes \rho_s
^{\ast}(\det(\fE)^{-1})$ is a quotient of $\rho_s^{\ast} (\fF)$,
it is nef; moreover pulling it back to $\PP Q_{s,\fE}$ it has a
surjection onto $\Oc_{\PP Q_{s,\fE}}(r) \otimes \pi_s^{\ast}(\det
(\fE))^{-1} \simeq \cO_{\PP Q_{s,\fE}}(r \lambda_{s, \fE})$, so
that $\lambda_{s,\fE}$ is nef.

Conversely, if $\lambda_{s,\fE}$ is nef, the $\Q$-Higgs bundle
$Q_{s,\fE} \otimes \rho_s ^{\ast}(\det(\fE)^{-1/r})$ is nef. Since
this is true for every $s$, the $\Q$-Higgs bundle $\fE \otimes
\det(\fE)^{-1/r}$ is H-nef; by taking the $r$-th symmetric power
we obtain that $\fF$ is H-nef.
\end{remark}

\begin{example} \label{exmagnum}
We give an example of an H-nef Higgs bundle which is not nef as an
ordinary bundle. Let $X$ be a projective surface of general type
that saturates Miyaoka-Yau's inequality, i.e., $3c_2(X) =
c_1(X)^2$ (surfaces of general type satisfying this condition are
exactly those that are uniformized by the unit ball in $\C^2$
\cite{S88}). The Higgs bundle $\fE$ whose underlying vector bundle
is $E = \Omega_X \oplus \Oc_X$ with the Higgs morphism
$\phi(\omega, f) = (0, \omega)$ is semistable and satisfies
$\Delta(E)=0$, so that the Higgs bundle $\fF=S^3(\fE)\otimes (\det
\fE)^{-1}$ is 1-H-nef. On the other hand,  the underlying vector
bundle $F = S^3(\Omega_X \oplus \Oc_X) \otimes K^{\ast}_X$
contains $K^{\ast}_X$ as a direct summand and therefore is not nef
(note that we exclude that $K_X$ is numerically flat).
\end{example}

From this Example we see that if $X$ is a projective surface of
general type for which  the Miyaoka-Yau inequality holds strictly,
i.e., $3c_2(X) > c_1(X)^2$, then there are curves $C$ in $X$ such
that the restriction of the Simpson system $E = \Omega_X \oplus
\Oc_X$ to $C$ is not semistable. It would interesting to investigate what
this tells us about the geometry of such curves.

It is interesting to note that the criterion expressed in Theorem
\ref{critproj} allows one to relate the semistability of a vector
bundle $E$ with the nefness of the tangent bundle to the total
space of the Grassmann varieties of $E$.

\begin{example}
Let $X$ be a complex projective manifold, and $E$ a rank $r$ vector
bundle on $X$. Let us choose an integer $s$ such that
$0<s<r$, and consider the exact sequence
$$0 \to T_{\grass_s(E)/X} \to T_{\grass_s(E)}
\to p_s^{\ast}(T_X) \to 0\,.$$
We have the following results:
\begin{enumerate} \item
If $E$ is semistable and $\Delta(E)=0$, and
$T_X$ is nef, then $ T_{\grass_s(E)}$ is nef. Since the tangent bundles to the fibres of the projection $T_{\grass_s(E)} \to T_X$ are nef, the  fact the $E$ is semistable and has vanishing discriminant  appear to be the conditions for these bundles to glue to a nef
bundle on $\grass_s(E)$.
\item If $ T_{\grass_s(E)}$ and $K_X^{-1}$ are  nef, then
$E$ is semistable and $\Delta(E)=0$.
\item As a consequence, when $T_X$ is numerically flat (e.g., $X$ is an Abelian variety,
or a hyperelliptic surface) then $ T_{\grass_s(E)}$ is nef
if and only if $E$ is semistable and $\Delta(E)=0$.
\end{enumerate}
\end{example}

\end{document}